%% file: init.tex
\documentclass{amsart}
\usepackage[utf8]{inputenc}

\input{shared-preamble.sty}
\usepackage[backend=biber]{biblatex}

\DeclareEmphSequence{\bfseries}

\author{Matt Alexander}
\date{\today}

\address{University of Regina\
  3737 Wascana Parkway\
  Regina, Saskatchewan, S4S 0A2\
  Canada}
\email{mpa097@uregina.ca}

\title{A Lie Correspondence in the Setting of Hopf and Frobenius Algebras}

\begin{document}

\begin{abstract}
  In this paper we demonstrate that Hopf algebras and Frobenius algebras provide an appropriate setting for an analogue of the classical Lie correspondence.
  We introduce a compatibility structure between Hopf and Frobenius algebras called Hopf-Frobenius modules,
  and give constructions that can be applied to any such module, which produce corresponding Lie groups and Lie algebras.
  Further, we show that classical Lie algebras and Lie groups give rise to Hopf-Frobenius modules for which the corresponding Lie groups and Lie algebras produced by our construction are those of the ordinary Lie correspondence.
\end{abstract}


\maketitle
\tableofcontents

\input{content}

\printbibliography

\end{document}

%% file: content.tex
\section{Introduction}
\label{sec:intro}

Lie algebras and Lie groups have played a central role in the study of Hopf algebras since Hopf's seminal paper in 1941 \cite{hopf-alg-first}. Hopf studied H-spaces (which he called \(  \Gamma   \)-Mannigfaltigkeiten), and in particular, Lie groups, by examining the map induced from the product \(  G \times G \to G   \), onto the cohomology (with coefficients in a field):
\begin{equation*}
  H^\bullet(G) \to H^\bullet(G) \otimes H^{\bullet}(G). 
\end{equation*}
This map gives the cohomology the structure of a coalgebra, and is furthermore compatible with the cup product. The algebraic structure and relations on cohomology were later axiomatized into the modern concept of Hopf algebra. For more information on the history of Hopf algebras, see \cite{history-of-hopf-alg} and \cite[Section~3.11]{Cartier-Patras-hopf-book}.

Even today, some of the fundamental examples of Hopf algebras come from Lie theory: the universal enveloping algebra of a Lie algebra \(  U[\mathfrak{g}]   \)  and the group algebra of a Lie group (or any group) \(  kG   \). As the group algebra only relies on the discrete structure of a group, a kind of topological analogue provides another example of a Hopf algebra: the convolution algebra $\mathcal{D}'_{1} (G)$ of all distributions on a Lie group supported at the identity \(  1 \in G   \)  (where the test functions are the smooth functions with compact support $C_{c}^{\infty}(G)$). This ends up being isomorphic as a Hopf algebra to the universal enveloping algebra of the corresponding Lie algebra \cite[Theorem~3.7.1]{Cartier-Patras-hopf-book}:
\begin{equation*}
 \mathcal{D}_1'(G) \cong U[\mathfrak{g}]
\end{equation*}
and thus only carries information on the identity component of the Lie group.

In this paper, we show that keeping track of additional topological data, in the form of a Frobenius algebra, allows us to recover in some cases the entire Lie group from $U[\mathfrak{g}]$, not merely the identity component. More precisely, we introduce a new kind of algebraic structure which we call a \emph{Hopf-Frobenius module}, which can be associated to any Lie group or Lie algebra: Let $H$ be a Hopf algebra with antipode $S$, and $(F, \gen{\cdot, \cdot})$ be a (suitably generalized) Frobenius algebra. A (left) Hopf-Frobenius module $(H,F)$ is a left $H$-module, $F$, such that the following two conditions are satisfied:
\begin{align*}
  \sum (-1)^{\abs{x} \abs{h_{(2)}}} \gen{h_{(1)} \cdot x, h_{(2)} \cdot y} &= \varepsilon_H(h) \gen{x,y},  \\
  \gen{h \cdot x, y} &= \gen{x, S(h) \cdot y},
\end{align*}
for all $h \in H, x,y \in F$. In this way, the antipode acts as a kind of adjoint operation on $H$.

Moreover, we prove that every Hopf-Frobenius module \textit{has} an associated topological group and Lie algebra, which, in the case that our Hopf-Frobenius module is a classical Lie group or Lie algebra (of dimension greater than $2$), reproduce the classical Lie correspondence. Our main result is the following (see \Cref{thm:lie-group-matches-topologically}):

\begin{theorem*}{}{}
  Let $G$ be one of the following Lie groups:
\begin{equation*}
 \{ O(n) \mid n > 2 \} \cup \{ U(n) \mid n \geq 2 \} \cup \{ Sp(\mathbb{R}^{2n}, \mathbb{R}) \},
\end{equation*}
where $Sp(\mathbb{R}^{2n}, \mathbb{R})$ is the real symplectic group (the group of invertible linear transformations on $\mathbb{R}^{2n}$ that preserve the symplectic form), and let $\rho: G \to GL(F)$ be its standard representation on the corresponding vector space $F$ with bilinear forms $\{ \gen{\cdot,\cdot}_i \}_{i \in I}$. In particular, the action of $G$ preserves the forms:
\begin{equation*}
\gen{gx,gy}_i = \gen{x,y}_i, \text{ for all } g \in G, x,y \in F, i \in I.
\end{equation*}
If $\mathfrak{g}$ is the Lie algebra of $G$ and $U[\mathfrak{g}]$ is its universal enveloping algebra then
\begin{equation*}
\cat{G}_F^{U[\mathfrak{g}]} \cong G
\end{equation*}
as topological groups, where \(  \cat{G}_F^{U[\mathfrak{g}]}   \) is the topological group associated to the Hopf-Frobenius module \(  U[\mathfrak{g}]   \).
\end{theorem*}

Note that Hopf algebras are strictly algebraic structures, but our result in the above theorem is an isomorphism of topological groups. This topological structure arises out of the second component in our definition of Hopf-Frobenius module: a Frobenius algebra on which the Hopf algebra acts.

Many classical Lie groups and algebras have canonical representations on vector spaces, on which they preserve the relevant bilinear forms. We will show that each of these vector spaces can be given a Frobenius algebra structure whose Frobenius form is the relevant bilinear form preserved by the corresponding Lie group. Frobenius algebras can be endowed with a topology induced from their bilinear form, which we can pull back onto the Hopf algebra. This will allow us to define a topology not only on the Hopf and Frobenius algebras, but also on the corresponding group of the Hopf-Frobenius module, turning it into a topological group.

From \cite[Corollary~3.41]{warner1971}, a topological group that is also a topological manifold has at most one compatible smooth structure. The results of the present paper then show that the structure of classical Lie groups can be distilled down to purely algebraic data (of Hopf and Frobenius algebras and appropriate compatibility conditions).

Furthermore, it is natural to ask how other Lie groups and Lie algebras fit into the Hopf-Frobenius module framework (exceptional Lie groups, Lie groups over finite or local fields, etc.), and what the corresponding topological groups and Lie algebras are for Hopf-Frobenius modules which don't come from universal enveloping algebras or group algebras. We will leave these computations to future work, but we add a small note at the end of the present paper about them.

\newthought{Structure of paper: } In \Cref{sec:prelim} we recall the standard notions from the theory of Hopf algebras, as well as establish the generalized version of Frobenius algebras that will be revelant for our Lie correspondence analogue. In \Cref{sec:hf-modules}, we define the category of Hopf-Frobenius modules and establish their basic structure. \Cref{sec:hf-mod-exs} presents a menagerie of examples of Hopf-Frobenius modules, demonstrating their appearance in algebra, functional analyis, and geometry. Finally in \Cref{sec:lie-correspondence}, we prove our version of the Lie correspondence which holds for Hopf-Frobenius modules, and show that it recovers the ordinary Lie correspondence for classical Lie groups and Lie algebras. In particular, we show that each Hopf-Frobenius module has associated Lie algebra (\Cref{lem:lie-alg-of-hopf}) and topological group (\Cref{prop:lie-group-of-hopf} and \Cref{prop:top-group}), and that this group forms a new Hopf-Frobenius module structure on our Frobenius algebra (\Cref{lem:lie-group-gives-hfm}), on which it acts continuously (\Cref{prop:action-cont}).

\newthought{Acknowledgements: } The author would like to thank Martin Frankland, Mart\'{i}n Argerami, Pierre Ouimet, Donald Stanley and Donald Yau for their advice and feedback. This work was supported by an NSERC Postgraduate Scholarship.

\section{Preliminaries}
\label{sec:prelim}

\newthought{Convention: Grading} The algebraic structures we examine below can be formulated in both the \textit{graded} and \textit{ungraded} settings. For the purposes of our Lie theory results, we will only need the ungraded case. Therefore, in what follows we shall make note of the graded formulation in our definitions, and use it when proving the basic theory of Hopf-Frobenius modules, but will specialize to the ungraded case once we consider the Lie theory in \Cref{sec:lie-correspondence}.

\subsection{Hopf algebras}

We begin by recalling the notion of Hopf algebra, before discussing Frobenius algebras and introducing our more general notion of them.

\newthought{Convention: Base field} The definitions of Frobenius and Hopf algebras we present can be formulated over \textit{general commutative rings}, but for the purposes of this paper, we will focus on the real and complex numbers. Thus throughout we will work over a field $k.$ 

\begin{definition}{Bialgebras and Hopf Algebras}{hopf-alg}
  A graded \textbf{bialgebra} $B$ is a graded $k$-algebra with linear maps $\Delta: B \rightarrow B \otimes B$ and $\varepsilon: B \rightarrow k$  (called the \textbf{comultiplication} (or coproduct) and \textbf{counit}, respectively) such that
  \begin{enumerate}
  \item $\Delta$ and $\varepsilon$ are graded algebra homorphisms of degree zero.
  \end{enumerate}
  A graded \textbf{Hopf algebra} is a bialgebra $H$ along with a degree zero graded linear map $S: H \to H$ (called the \textbf{antipode}) such that
  \begin{equation*}
   \eta \circ \varepsilon(x) = \sum x_{(1)} S(x_{(2)}) = \sum S(x_{(1)}) x_{(2)}.
  \end{equation*}

\end{definition}
\begin{defframe}

\newthought{Involution: } If $S$ is an involution (i.e. $S^{2}= id)$, we call $H$ an \emph{involutive Hopf algebra}. In what follows all of our Hopf algebras will be assumed involutive unless otherwise specified.
\end{defframe}

\begin{definition}{Grouplike Elements and Primitive Elements}{grouplike-primitive}
Let $B$ be a bialgebra. An element $g \in B$ is called \emph{grouplike} if it has the comultiplication
\begin{align*}
\Delta g = g \otimes g, \quad g \ne 0.
\end{align*}
An element $p \in B$ is called \emph{primitive} if
\begin{align*}
\Delta p = 1 \otimes p + p \otimes 1.
\end{align*}
\end{definition}

\begin{proposition}{Properties of Grouplike and Primitive Elements}{•}
Let $B$ be a bialgebra.
\begin{enumerate}
\item The grouplike elements of $B$ form a monoid. Furthermore, if $B$ is a Hopf algebra with antipode $S$, the monoid of grouplike elements is a group, where the inverse of $g$ is $S(g)$.

\item The primitive elements of $B$ form a Lie algebra under the commutator bracket. 
\end{enumerate}
\end{proposition}
\begin{proof-prop}
See \cite[Proposition~5.1.15]{radfordtextbook}.
\end{proof-prop}

  \newthought{Examples: } The following two examples show how Hopf algebras can be built out of grouplike and primitive elements.

\begin{example}{Hopf Algebra: Group Algebra}{hopf-group-alg}

  Let $G$ be a group (finite or infinite) and $k$ be a field. The group algebra $kG$ has the structure of a Hopf algebra with
  \begin{equation}
\begin{split}
\Delta g &= g \otimes g \\
\varepsilon(g) &= 1 \\
S(g) &= g^{-1}.
\end{split}
\end{equation}
for all $x \in G$.
\end{example}

\begin{example}{Hopf Algebra: Universal Enveloping Algebra}{hopf-universal-envel}
Let $\mathfrak{g}$ be a Lie algebra. The universal enveloping algebra $U(\mathfrak{g})$ has the structure of a Hopf algebra with
\begin{equation}
\begin{split}
\Delta x &= x \otimes 1 + 1 \otimes x \\
\varepsilon(x) &= 0 \\
S(x) &= -x.
\end{split}
\end{equation}
for all $x \in \mathfrak{g}$.
\end{example}

\subsection{Frobenius algebras}

Frobenius algebras are ordinarily taken to be strictly associative and finite-dimensional. For our purposes we will consider a weaker notion of Frobenius algebra, which relaxes associativity as well as permits an extra involution on the algebra. These features will allow us to treat Lie groups involving sesquilinear, orthogonal, symplectic, and exceptional structures all within the same framework.

\begin{definition}{$*$-ring}{}
A \emph{$*$-ring} is a ring $R$ with an involution $*: R \to R$ such that
\begin{align*}
(x + y)^* &= x^* + y^* \\
(xy)^* &= y^* x^* \\
1^* &= 1,
\end{align*}
for all $x, y \in R$.
\end{definition}
\begin{defframe}

\newthought{Note: } Every commutative ring is a $*$-ring with trivial involution.

\newthought{Convention: } In what follows we will consider $k$ to be a $*$-field unless otherwise specified.
\end{defframe}

\begin{definition}{Graded $*$-Algebra}{}
A \emph{graded $*$-algebra} is a graded $k$-algebra $A$ with an involution \\ $*: A \to A$ such that 
\begin{align*}
(\alpha x + \beta y)^* &= \alpha^* x^* + \beta^* y^* \\
(xy)^* &= x^* y^*,
\end{align*}
for all $x,y \in A$ and $\alpha, \beta \in k$.
\end{definition}
\begin{defframe}

\newthought{Convention: } The involution in a $*$-algebra is often required to satisfy $(xy)^* = y^* x^*$. However, in order to treat ordinary algebras (over ordinary, non-starred rings) and $*$-algebras on the same footing, we will replace this requirement by multiplicativity of the involution.
\end{defframe}

\begin{example}{Graded $*$-Algebra: Graded Algebra}{}
Every algebra $A$ over an ordinary (trivially starred) commutative ring has a trivial involution $x^* = x$ for all $x \in A$. Thus we can view every algebra as a $*$-algebra according to our definition.
\end{example}

\begin{example}{Graded $*$-algebra: Complex Matrices with Entrywise Conjugation}{}

Let $M_{n}(\mathbb{C})$ be the algebra of $n \times n$ complex matrices, and take as involution the entrywise conjugation of matrices. This defines a $*$-algebra in our sense. 

\end{example}
\begin{exampleframe}
We will refer back to this example in \Cref{ex:frb-matrices-with-conj}, where we show that the Frobenius structure on this algebra reproduces the ordinary inner product on complex matrices. 
\end{exampleframe}

\begin{definition}{Graded Frobenius Algebra}{}
A graded \emph{Frobenius $*$-algebra} $(F, \gen{\cdot, \cdot})$, is a graded not necessarily associative algebra $F$ over a $*$-ring $k$, with non-degenerate bilinear form $\gen{\cdot, \cdot}: F \otimes F \to k$, such that
\begin{equation}
\gen{xy,z} = \gen{x,y z}
\end{equation}
for all $x,y,z \in F$, and such that for all non-zero $x \in F$, there exist some $a,b \in F$ with $a x \ne 0$ and $x b \ne 0$. We do not assume that F is unital nor finite-dimensional.
\end{definition}
\begin{defframe}
Our notion of Frobenius algebra generalizes the standard notion in four ways, permitting: non-associative algebras, infinite-dimensional algebras, a graded structure, and an involution.

\newthought{Trivial Structure Case: } In the case of a finite-dimensional associative algebra with trivial involution, we recover the usual notion of Frobenius algebra.

\newthought{Non-Associativity: } The condition on the pairing, $\gen{xy,z} = \gen{x,y z}$ can be thought of as a kind of weak associativity condition. We will see in \Cref{lem:frob-unital-counit} that under certain assumptions, this weak associativity condition is equivalent to the existence of a linear functional $\varepsilon$ on the algebra such that $\varepsilon((xy)z) = \varepsilon(x(yz))$ for all $x,y,z \in F$.

\newthought{Infinite Dimensions: } In the literature, infinite-dimensional Frobenius algebras are often defined differently to the approach we take here, requiring a comultiplication $\Delta: F \to F \otimes F$, required to satisfy certain compatibility conditions with the multiplication on $F$. For more information, see \cite{frob-nearly, frob-non-unital}.

\newthought{Involution: } Defining Frobenius algebras and involutions in the way we have will allow us to treat real and sesquilinear pairings as two sides of the same structure. We make this point more explicitly in \Cref{lem:frob-2-pairings} below.
\end{defframe}

\medskip

We first show that the structure of a Frobenius algebra can be expressed in terms of either a bilinear or sesquilinear form, satisfying the appropriate associativity conditions. 

\begin{lemma}{}{frob-2-pairings}
Let $F$ be a graded Frobenius $*$-algebra with bilinear form $\gen{\cdot, \cdot}_{\text{Frob}}$. Then there is an associated sesquilinear form $\gen{\cdot, \cdot}$, defined by 
\begin{equation}
\gen{x,y} \equiv \gen{x,y^*}_{\text{Frob}}
\end{equation}
which satisfies
\begin{equation}\label{eq:sesq-assoc}
\gen{xy, z} = \gen{x,y^* z}.
\end{equation}
Similarly, given a sesequilinear form on $F$, $\gen{\cdot, \cdot}$, satisfying \Cref{eq:sesq-assoc}, there is an associated bilinear form $\gen{\cdot, \cdot}_{\text{Frob}}$ that gives $F$ the structure of a Frobenius algebra, defined by
\begin{equation}
\gen{x,y}_{\text{Frob}} \equiv \gen{x, y^*}.
\end{equation}
\end{lemma}
\begin{proof-lemma}

\newthought{If Bilinear: }  Given a Frobenius form $\gen{\cdot, \cdot}_{\text{Frob}}$, we have
\begin{align*}
\gen{xy, z} &= \gen{xy, z^*}_{\text{Frob}} \\
&= \gen{x, y z^*}_{\text{Frob}} \\
&= \gen{x, (y^*)^* z^*}_{\text{Frob}} \\
&= \gen{x, (y^* z)^*}_{\text{Frob}}  \\
&= \gen{x, y^* z}.
\end{align*}
The induced pairing is sesquilinear, as 
\begin{align*}
\gen{\alpha x + \beta y, \gamma z + \delta w} &= \gen{\alpha x + \beta y, (\gamma z + \delta w)^*}_{\text{Frob}} \\
&= \gen{\alpha x + \beta y, \gamma^* z^* + \delta^* w^*}_{\text{Frob}} \\
&= \alpha \gamma^* \gen{x,z^*}_{\text{Frob}} + \alpha \delta^* \gen{x,w^*}_{\text{Frob}} + \beta \gamma^* \gen{y, z^*}_{\text{Frob}} + \beta \delta^* \gen{y, w^*}_{\text{Frob}} \\
&= \alpha \gamma^* \gen{x,z} + \alpha \delta^* \gen{x,w} + \beta \gamma^* \gen{y, z} + \beta \delta^* \gen{y, w}.
\end{align*}

\newthought{If sesquilinear: } A symmetric argument produces a Frobenius form from a sesquilinear one satisfying \Cref{eq:sesq-assoc}.
\end{proof-lemma}

\medskip

Often a Frobenius algebra can be defined purely in terms of a linear functional $\varepsilon : F \rightarrow k$, which we call a \textbf{counit}. This will come in handy for the examples in \Cref{sec:frb-ex}

\begin{lemma}{}{frob-unital-counit}
Let $F$ be an \textbf{associative} unital graded Frobenius algebra. Then there is a linear functional $\varepsilon: F \to k$ whose kernel contains no non-trivial ideals (left nor right), such that 
\begin{equation}
\gen{x,y} = \varepsilon(xy^*).
\end{equation}
Conversely, given any non-zero functional $\varepsilon$ with simple kernel on an associative algebra, $\gen{x,y} = \varepsilon(xy^*)$ defines a (sesquilinear) Frobenius form.
\end{lemma}
\begin{proof-lemma}

\newthought{Given a form: } Since $F$ is unital, we can define a functional $\varepsilon(y) \equiv \gen{1,y^*}$. Then $\varepsilon(xy^*) = \gen{1,(xy^*)^*} = \gen{1,x^* y} = \gen{x,y}$. Now let $I$ be a right ideal in the kernel of $\varepsilon$ and let $x \in I$. It follows that $xy^* \in I$ for all $y$, so $\gen{x,y} = \varepsilon(xy^*) = 0$ for all $y$. Thus $x$ must be zero, and $I$ trivial. Similarly, if $I$ is a left ideal, taking $y \in I$ would imply $\gen{x,y^*} = \varepsilon(xy) = 0$ for all $x$, and so $I$ must be trivial.

\newthought{Given a functional: } $\gen{xy,z} = \varepsilon(xy z^*) = \varepsilon(x(y^* z)^*) = \gen{x,y^* z}$. Now if $\gen{x,y} = 0$ for all $x$, we can take the left ideal generated by $y^*, I_{y^*}^L$, which will be contained in the kernel of $\varepsilon$. Since $\varepsilon$ is non-zero by assumption, it follows that $I_{y^*}^L$ must be trivial, so $y^* = 0$. Using the linearity of the involution, we conclude $y = y^{**} = 0$. Similarly, if $\gen{x,y} = 0$ for all $y$, the right ideal generated by $x$ must be trivial, and thus $x$ must be zero. It follows that the pairing is non-degenerate.
\end{proof-lemma}

\subsection{Examples of Frobenius algebras}
\label{sec:frb-ex}

For the purposes of this paper, there are three main examples of Frobenius algebras that we need to consider. Each is constructed in a similar way: by picking a basis on \(  \mathbb{R}^n   \) and defining an appropriate multiplication on the basis elements, which recovers the desired bilinear form once the counit is applied. Before constructing those examples, however, we look at two examples which will motivate some of the assumptions we imposed in our definition of Frobenius algebras.

\Cref{ex:frb-matrices-with-conj} motivates taking the involution on our algebras to be linear, rather than anti-linear, and \Cref{ex:comp-alg} (particularly the case of the octonions) motivates allowing non-associative algebras.

\begin{example}{Frobenius $*$-Algebra: $n \times n$ Matrices}{frb-matrices-with-conj}

  Let $(M_{n}(\mathbb{C}), *)$ be the $*$-algebra of $n \times n$ matrices with entrywise conjugation. We define our bilinear form as
  \begin{equation*}
   \gen{A,B}_{\text{Frob}} = \operatorname{Tr}(A B^{T}).
  \end{equation*}
The corresponding sesquilinear pairing is then the ordinary inner product of matrices:
\begin{equation*}
 \gen{A,B} = \operatorname{Tr}(A B^{\dagger}).
\end{equation*}
\end{example}

\begin{example}{Frobenius $*$-Algebra: Composition Algebras}{comp-alg}
A \emph{composition algebra} is a finite-dimensional vector space $V$ with not necessarily associative multiplication, a unit, and a non-degenerate symmetric bilinear form $(\cdot, \cdot)$, such that 
\begin{equation}
(xy, xy) = (x,x) (y,y).
\end{equation}

\end{example}
\begin{exampleframe}

  \newthought{Convention: } A composition algebra in the sense defined above yields a composition algebra in the sense of \cite[Chapter~1]{springer-octonions} with quadratic form $N(x) = (x,x)$, where they define a bilinear form $B(x,y) = N(x+y) - N(x) - N(y)$.  Assuming $char(k) \neq 2$, our bilinear form will satisfy
  \begin{equation*}
 2 (x,y) = B(x,y). 
  \end{equation*}
and the two notions are equivalent.

\newthought{Composition involution: } Every composition algebra has an involution defined by
\begin{equation}
x^\dagger \equiv 2 (x,1)1 - x.
\end{equation}
It turns out that $(x y, z) = (x,zy^\dagger)$ and $(xy)^\dagger = y^\dagger x^\dagger$. We can define a new non-degenerate pairing on any composition algebra by
\begin{equation}
\gen{x,y}_{\text{Frob}} \equiv (x, y^\dagger),
\end{equation}
which will satisfy $\gen{xy,z}_{\text{Frob}} = \gen{x, yz}_{\text{Frob}} $.

\newthought{Frobenius Involution: } Note that the involution on the composition algebra is anti-multiplicative: $(xy)^\dagger = y^\dagger x^\dagger$, but the involutions on our Frobenius algebras are defined to be \textit{multiplicative}. So we view our underlying ring $k$ as a $*$-ring with trivial involution and take our involution on the Frobenius algebra to be trivial as well.

\newthought{Complex Numbers: } The complex numbers can be viewed as a real composition algebra, where 
\begin{align*}
(x,y) \equiv \dfrac{1}{2} (x^* y + xy^*).
\end{align*}
The involution on this algebra is the complex conjugate: $x^\dagger = x^*$, and the Frobenius form we end up with is the real part of the product:
\begin{align*}
\gen{x,y}_{\text{Frob}} = \mathbb{R}(xy).
\end{align*}
On the other hand, taking the complex numbers as a complex composition algebra, where
\begin{align*}
(x,y) = xy
\end{align*}
and the trivial involution is $x^\dagger = x$ gives
\begin{align*}
\gen{x,y}_{\text{Frob}} = xy.
\end{align*}
Note in this case, even though the composition algebra involution $(-)^\dagger$ is trivial, on our complex Frobenius algebra we can take the complex conjugate as our involution, giving
\begin{align*}
\gen{x,y} = xy^*,
\end{align*}
the ordinary complex inner product. Thus the involutions in the composition and Frobenius algebras play different roles in general.

\newthought{Octonions: } The octonions can be viewed as a real composition algebra as follows: let $\{ e_0, \ldots, e_7 \}$ be the standard basis of the octonions, and write each $x = \sum_i x_i e_i$. Then
\begin{align*}
(x,y) = \sum_i x_i y_i
\end{align*}
and $x^\dagger = x_0 e_0 - \sum_{i=1}^7 x_i e_i$. We get
\begin{align*}
\gen{x,y}_{\text{Frob}} = x_0 y_0 - \sum_{i=1}^7 x_i y_i.
\end{align*}
This provides a nice example of a Frobenius algebra which is not associative, but still satisfies associativity under the pairing:
\begin{align*}
\gen{xy, z} = \gen{x, yz}.
\end{align*}
We can think of Frobenius algebras of this type as being \textit{weakly associative}.

\end{exampleframe}

\begin{example}{Frobenius Algebra: Euclidean Space}{frb-euclid}

  Let $V$ be a finite-dimensional real vector space and choose a basis $\{ e_i \}$. We can give $V$ the structure of a Frobenius algebra with multiplication
  \begin{equation*}
   e_i e_j = \delta_{ij}e_j
  \end{equation*}
  and pairing
  \begin{equation*}
   \gen{e_i, e_j} = \delta_{ij}
  \end{equation*}
  where $\delta_{ij}$ is the Kronecker delta.
\end{example}
\begin{exampleframe}
  A pairing of two general elements $x = \sum_i x_i e_i$ and $y = \sum_i y_i e_i$ is the ordinary pairing defined by the basis:
  \begin{equation*}
   \gen{x,y} = \sum_{i}x_i y_i.
  \end{equation*}
  \newthought{Generalization: } We can replace the base field by $\mathbb{C}$ and use complex conjugation as our involution to acquire the ordinary sesquilinear pairing:
  \begin{equation*}
    \gen{x,y} = \gen{x,y^{*}}_{\text{Frob}} = \sum_{i} x_{i} y_{i}^{*}.
  \end{equation*}
  By the same token, replacing $\mathbb{R}$ by the quaternions $\mathbb{H}$ with conjugation, produces a Frobenius $*$-algebra with the ordinary quaternionic pairing.
\end{exampleframe}

\medskip

The same underlying vector space can permit a number of different Frobenius algebra structures. Consider the following three examples.

\begin{example}{Real Frobenius $*$-Algebra: Complex Numbers}{frb-complex}
Consider the real vector space $\mathbb{R}^{2n}$, with basis $\{ a_i, b_i \}_{i=1}^n$. We impose the following multiplication
\begin{align*}
a_i b_j &= \delta_{ij} b_i \\
b_i a_j &= \delta_{ij} b_i \\
a_i a_j &= \delta_{ij} a_i \\
b_i b_j &= -\delta_{ij} a_i,
\end{align*}
take the counit to be $\varepsilon(a_i) = 1$, $\varepsilon(b_i) = 1$, and let the star involution be $a_i^* = a_i$, $b_i^* = - b_i$.
\end{example}
\begin{exampleframe}
For any fixed $i$, we have
\begin{align*}
(\alpha_1 a_i + \beta_1 b_i)(\alpha_2 a_i + \beta_2 b_i) = (\alpha_1 \alpha_2 - \beta_1 \beta_2) a_i + (\alpha_1 \beta_2 + \beta_1 \alpha_2) b_i,
\end{align*}
recovering the complex numbers. The induced Frobenius pairing will be
\begin{align*}
\gen{\alpha_i a_i + \beta_i b_i, (\alpha_j a_j + \beta_j b_j)^*}_{\text{Frob}} &=  \delta_{ij}(\alpha_i \alpha_j + \beta_i \beta_j - \alpha_i \beta_j + \beta_i \alpha_j).
\end{align*}
If we view $z_i = \alpha_i a_i + \beta_i b_i$ as an element of $\mathbb{C}^n$, the above pairing is the composition of the ordinary pairing on $\mathbb{C}^n$ with the $\mathbb{R}$-linear map $\pi: \mathbb{C} \to \mathbb{R}$ that sends $\pi(1) = \pi(i) = 1$. 
\end{exampleframe}

\begin{example}{Real Frobenius $*$-Algebra: Orthogonal Form}{frb-orthog}
Take the same $*$-algebra structure on $\mathbb{R}^{2n}$ as in the previous example, but consider the following counit:
\begin{align*}
\varepsilon_O(a_i) &= 1 \\
\varepsilon_O(b_i) &= 0.
\end{align*}
The Frobenius form in this example is the standard orthogonal form on $\mathbb{R}^{2n}$.
\end{example}
\begin{exampleframe}
From our previous example, we know that 
\begin{align*}
(\alpha_1 a_i + \beta_1 b_i)(\alpha_2 a_i + \beta_2 b_i) = (\alpha_1 \alpha_2 - \beta_1 \beta_2) a_i + (\alpha_1 \beta_2 + \beta_1 \alpha_2) b_i.
\end{align*}
It follows that $\varepsilon((\alpha_1 a_i + \beta_1 b_i)(\alpha_2 a_i + \beta_2 b_i)^*) = \alpha_1 \alpha_2 + \beta_1 \beta_2$. Thus
\begin{align*}
\varepsilon_O(xy^*) = \gen{x,y}_O = \sum_i x_i y_i,
\end{align*}
yielding the standard inner product.
\end{exampleframe}

\begin{example}{Real Frobenius $*$-Algebra: Symplectic Form}{frb-symp}
If we take the same $*$-algebra structure on $\mathbb{R}^{2n}$ as in the previous two examples, but select our counit to be 
\begin{align*}
\varepsilon_S(a_i) &= 0 \\
\varepsilon_S(b_i) &= 1
\end{align*}
our Frobenius form becomes the standard symplectic form.
\end{example}
\begin{exampleframe}
The Frobenius form $\gen{x,y} = \varepsilon(x y^*)$ satisfies
\begin{align*}
\gen{a_i, a_j} &= 0  & \gen{a_i, b_j} &= - \delta_{ij} \\
\gen{b_i, b_j} &= 0 & \gen{b_i, a_j} &= \delta_{ij},
\end{align*}
giving the standard symplectic pairing.
\end{exampleframe}

\newthought{Unitary matrices: } The linear maps which simultaneously preserve both the Frobenius structures on $\mathbb{R}^{2n}$ in \Cref{ex:frb-orthog} and \Cref{ex:frb-symp} can be identified with the unitary matrices $U(n)$. This result is part of the \textbf{2-out-of-3 property} (see for instance \cite[Section 41]{arnold-2outof3}). This suggests that considering multiple Frobenius algebra structures on the same underlying algebra can be useful. We will utilize this perspective in \Cref{def:lie-group-hf-mod}

\subsection{Examples from functional analysis}

To motivate the infinite-dimensional generalization of Frobenius algebras, we make note of several examples from functional analysis.

\begin{example}{Frobenius Algebra: Hilbert Space with Basis}{frb-hilbert}
Let $H$ be a Hilbert space (even one that is non-separable). Pick an orthonormal basis $\{ e_\alpha \}$ and define a multiplication as follows: given any two elements $x = \sum_{i=1}^{\infty} x_{\alpha_i} e_{\alpha_i}$, $y = \sum_{j=1}^{\infty} y_{\beta_j} e_{\beta_j}$, let 
\begin{align*}
xy \equiv \sum_{k=1}^{\infty} x_{\gamma_k} y_{\gamma_k} e_{\gamma_k},
\end{align*}
where $\{ e_{\gamma_k} \}_{k=1}^{\infty} = \{ e_{\alpha_i} \}_{i=1}^{\infty} \cup \{ e_{\beta_j} \}_{j=1}^{\infty}$, and the coefficients $x_i, y_j$ have been reindexed appropriately. With the pairing
\begin{align*}
\gen{x,y}_{\text{Frob}} \equiv \sum_{k=1}^{\infty} x_{\gamma_k} y_{\gamma_k},
\end{align*}
this has the structure of a Frobenius algebra. For more information on non-separable Hilbert spaces, see \cite{martin-book}.
\end{example}
\begin{exampleframe}
We need to check that the multiplication is well-defined: that it converges and doesn't depend on how we chose to order the $e_{\gamma_k}$ in the sum. 

\newthought{Pythagorean Theorem: } Since $\{ e_\alpha \}$ are orthonormal, we have
\begin{align*}
\norm{\sum_{i=1}^N \lambda_i e_{\alpha_i}}^2 &=  \gen{\sum_{i=1}^N \lambda_i e_{\alpha_i}, \sum_{i=1}^N \lambda_i e_{\alpha_i} }_H = \sum_{i,j=1}^N \gen{\lambda_i e_{\alpha_i}, \lambda_j e_{\alpha_j}}_H \\
&= \sum_{i=1}^N \abs{\lambda_i}^2. 
\end{align*}

\newthought{Convergence: } First note if $x = \sum_i x_i e_{\alpha_i}$, then $( \sum_{i=1}^N x_i e_{\alpha_i} )_{N=1}^{\infty}$ is a Cauchy sequence. Considering the differences of the partial sums, we have
\begin{align*}
\norm{\sum_{i=1}^N x_i e_{\alpha_i} - \sum_{i=1}^M x_i e_{\alpha_i}} &= \norm{\sum_{i=M+1}^N x_i e_{\alpha_i}} \\
&= \sum_{i=M+1}^N \abs{x_i}^2.
\end{align*}
It follows that the sequence $( \sum_{i=1}^N x_i e_{\alpha_i} )$ is Cauchy if and only if $(\sum_{i=1}^N \abs{x_i}^2)$ is Cauchy. Now if $x$ and $y$ are expressed in terms of our orthonormal basis, we can show that $\sum_k x_{\gamma_k} y_{\gamma_k} e_{\gamma_k}$ is absolutely convergent: for any fixed $M,N$, by the Cauchy-Schwarz inequality,
\begin{align*}
\sum_{k=M}^N \abs{x_{\gamma_k} y_{\gamma_k}} \leq \sqrt{\sum_{k=M}^N \abs{x_{\gamma_k}}^2 \sum_{\ell=M}^N \abs{y_{\gamma_\ell}}^2  }.
\end{align*}
Since $(\sum_{k=1}^N \abs{x_{\gamma_k}}^2)$ and $(\sum_{\ell=1}^N \abs{y_{\gamma_\ell}}^2)$ are Cauchy by the argument above, we can make the right hand side of the inequality arbitrarily small by picking $M,N$ sufficiently large. Thus the series $\sum_{k} \abs{x_{\gamma_k} y_{\gamma_k}}$ is absolutely convergent. In a Hilbert space, absolute convergence implies unconditional convergence. Thus the series
\begin{align*}
\sum_{k=1}^{\infty} x_{\gamma_k} y_{\gamma_k} e_{\gamma_k}
\end{align*}
is both convergent and independent of how the terms in the sum are rearranged. So our multiplication is well-defined.

\newthought{Non-degeneracy: } First note that our pairing 
\begin{align*}
\gen{x,y}_{\text{Frob}} \equiv \sum_{k=1}^{\infty} x_{\gamma_k} y_{\gamma_k}
\end{align*}
is well-defined, since the sum on the right is absolutely convergent from the argument above, and thus converges. Note that the Frobenius pairing matches the inner product on the Hilbert space. Now given any non-zero $x$, there is some coefficient $x_{\alpha_j} \ne 0$. Then
\begin{align*}
\gen{x, e_{\alpha_j}} = x_{\alpha_j} \ne 0,
\end{align*}
and so our pairing is non-degenerate.

\newthought{Associativity: } 
\begin{align*}
\gen{x,yz}_{\text{Frob}} &= \gen{x, \sum_k y_{\gamma_k} z_{\gamma_k} e_{\gamma_k}} = \sum_\ell x_{\theta_{\ell}} y_{\theta_{\ell}} z_{\theta_{\ell}}\\
&= \gen{xy,z}_{\text{Frob}}.
\end{align*}
Both pairings are equal to the sum of the product of those coefficients of $x,y,z$ that are non-zero on the same basis elements.

\end{exampleframe}

\begin{example}{Frobenius Algebra: Compactly Supported Functions}{frb-comp-supp}
Let $C_c^{\infty}(\mathbb{R}^n)$ be the algebra of smooth compactly supported functions. The pairing
\begin{equation}
\gen{f,g}_{\text{Frob}} = \int f g
\end{equation}
gives $C_c^{\infty}(\mathbb{R}^n)$ the structure of a Frobenius algebra.
\end{example}
\begin{exampleframe}
This example works more generally: Let $M$ be any oriented smooth manifold equipped with a volume form $\omega \in \Omega^n(M)$. The chosen volume form gives us a way of defining an integral over the manifold. Taking the algebra $C_c^{\infty}(M)$ of smooth compactly supported functions, we can define a Frobenius form by
\begin{align*}
\gen{f,g}_{\text{Frob}} = \int_M f g \omega,
\end{align*}
giving $C_c^{\infty}(M)$ the structure of a Frobenius algebra.
\end{exampleframe}

\begin{example}{Frobenius Algebra: Schwartz Space}{frb-schwartz}
Let 
\begin{align*}
\cat{S} = \{ f \in C^{\infty}(\mathbb{R}^n, \mathbb{C}) \mid \forall \alpha, \beta, \sup_{x \in \mathbb{R}^n} \abs{x_1^{\alpha_1} \ldots x_n^{\alpha_n} (\partial_1^{\beta_1} \ldots \partial_n^{\beta_n} f)} < \infty \}
\end{align*}
be the algebra of Schwartz functions. The pairing
\begin{equation}
\gen{f,g} = \int f g^*
\end{equation}
gives $\cat{S}$ the structure of a Frobenius $*$-algebra.
\end{example}
\begin{exampleframe}
  For an introduction to the theory of Schwartz spaces and distributions, see for example \cite[Chapter 10]{first-course-sobolev-spaces}.
\end{exampleframe}

\section{Hopf-Frobenius modules}
\label{sec:hf-modules}

\begin{definition}{Hopf-Frobenius Module}{hf-mod}
A \emph{left Hopf-Frobenius module} is a pair $(H,F)$ where $H$ is a (graded) Hopf algebra, and $F$ is a (graded) Frobenius algebra such that
\begin{enumerate}
\item $F$ has the structure of a left $H$-module.

\item $\sum (-1)^{\abs{x} \abs{g_{(2)}}} \gen{g_{(1)} x, g_{(2)} y} = \varepsilon_H(g) \gen{x,y}$.

\item $\gen{g x, y} = (-1)^{\abs{g} \abs{x}} \gen{x, S(g) y} $ 
\end{enumerate}
for all $g \in H$ and $x,y \in F$.
\end{definition}
\begin{defframe}

  \medskip

\newthought{Convention: } One can similarly define right Hopf-Frobenius modules. However, in what follows, we shall consider all of our Hopf-Frobenius modules to be \textit{left} Hopf-Frobenius modules.
\end{defframe}

The third condition in the definition of Hopf-Frobenius modules will be useful in proving a version of the Lie group -- Lie algebra correspondence for Hopf-Frobenius modules. However, for grouplike and primitive elements of the Hopf algebra, it follows from the previous conditions:

\begin{lemma}{}{}
Let $H$ be a Hopf algebra and $F$ be Frobenius algebras which is also a left $H$-module. If $H$ and $F$ satisfy
\begin{align*}
\sum (-1)^{\abs{x} \abs{g_{(2)}}} \gen{g_{(1)} x, g_{(2)} y} &= \varepsilon_H(g) \gen{x,y}
\end{align*}
and $g$ is grouplike or primitive, then $\gen{S(g) x,y} = (-1)^{\abs{g} \abs{x}} \gen{x,gy}$.
\end{lemma}
\begin{proof-lemma}

\newthought{Primitive: } Let $p$ be primitive. The first condition for HF-modules becomes 
\begin{align*}
  \sum (-1)^{\abs{x} \abs{p_{(2)}}} \gen{p_{(1)} x, p_{(2)} y} &= \gen{p x, y} + (-1)^{\abs{x} \abs{p}} \gen{x, py} \\
                                                               &= \varepsilon(p) \gen{x,y} \\
                                                               &= 0.
\end{align*}
It follows that
\begin{align*}
  \gen{px, y} &= - (-1)^{\abs{x} \abs{p}} \gen{x, py} \\
  &= (-1)^{\abs{x} \abs{p}} \gen{x, (-p) y} \\ &= (-1)^{\abs{x} \abs{p}} \gen{x, S(p) y}.
\end{align*}

\newthought{Grouplike: } Let $g$ be grouplike. The first condition becomes $(-1)^{\abs{x} \abs{g}} \gen{gx, gy} = \gen{x,y}$. Since $g^{-1}$ is also grouplike, we have $\gen{g x, y} = (-1)^{\abs{x} \abs{g^{-1}}} \gen{g^{-1} g x, g^{-1} y} = (-1)^{\abs{x} \abs{g}} \gen{x, S(g) y}$.
\end{proof-lemma}

\medskip

Since the relevant Hopf algebras for our purposes are generated linearly and multiplicatively through grouplike and primitive elements, the following lemma shows that satisfying condition (2) in \Cref{def:hf-mod} is the central criterion for our Hopf-Frobenius modules.

\begin{lemma}{}{}
Let $H$ be a Hopf algebra and $F$ be a Frobenius algebra that is also a left $H$-module. If $g,h \in H$ satisfy 
\begin{align*}
\sum (-1)^{\abs{x} \abs{g_{(2)}}} \gen{g_{(1)} x, g_{(2)} y} = \varepsilon_H(g) \gen{x,y}
\end{align*}
for all $x,y \in F$ then so does the product $gh$. Similarly, if $g,h$ satisfy
\begin{align*}
\gen{g x, y} = (-1)^{\abs{g} \abs{x}} \gen{x, S(g) y}
\end{align*}
for all $x,y$, 
then so does $gh$.
\end{lemma}
\begin{proof-lemma}

\newthought{Condition 2: } First note that
\begin{align*}
\Delta(gh) &= \Delta(g) \Delta(h) \\
&= \left( \sum g_{(1)} \otimes g_{(2)} \right) \left( \sum h_{(1')} \otimes h_{(2')} \right) \\
&= \sum (-1)^{\abs{g_{(2)}} \abs{h_{(1')}}} g_{(1)} h_{(1')} \otimes g_{(2)} h_{(2')}.
\end{align*}
Using this, we can write the second Hopf-Frobenius condition in \Cref{def:hf-mod} applied to the product $gh$ as
\begin{align*}
\sum (-1)^{\abs{x} \abs{(gh)_{(2)}}} &\gen{(gh)_{(1)} x, (gh)_{(2)} y} \\
&= \sum (-1)^{\abs{x} \abs{g_{(2)} h_{(2')}} + \abs{g_{(2)}} \abs{h_{(1')}}} \gen{(g_{(1)} h_{(1')}) x, (g_{(2)} h_{(2')}) y} \\
&= \sum (-1)^{\abs{x} \abs{g_{(2)}} +  \abs{x} \abs{ h_{(2')}} + \abs{g_{(2)}} \abs{h_{(1')}}} \gen{g_{(1)} (h_{(1')} x), g_{(2)} (h_{(2')} y)}\\
&= \sum (-1)^{ \abs{g_{(2)}}\abs{h_{(1')} x} +  \abs{x} \abs{ h_{(2')}} } \gen{g_{(1)} (h_{(1')} x), g_{(2)} (h_{(2')} y)} \\
&= \varepsilon(g) \sum (-1)^{  \abs{x} \abs{ h_{(2')}} } \gen{h_{(1')} x, h_{(2')} y} \\
&= \varepsilon(g) \varepsilon(h) \gen{x,y}\\
&= \varepsilon(gh) \gen{x,y}.
\end{align*}

\newthought{Condition 3: } Finally we check the third condition in \Cref{def:hf-mod}. We use the fact that 
\begin{align*}
S(gh) = (-1)^{\abs{g} \abs{h}} S(h) S(g).
\end{align*}
Then we have
\begin{align*}
\gen{(gh)x, y} &= \gen{g(hx), y}\\
&= (-1)^{\abs{g} \abs{hx}} \gen{hx, S(g) y}\\
&= (-1)^{\abs{g} \abs{h} + \abs{g} \abs{x} + \abs{h} \abs{x}} \gen{x, S(h)S(g) y} \\
&= (-1)^{\abs{g} \abs{x} + \abs{h} \abs{x}} \gen{x, S(gh) y} \\
&= (-1)^{\abs{x}\abs{gh}} \gen{x, S(gh) y}. \qedhere
\end{align*}
\end{proof-lemma}

\subsection{The inherited topology on Frobenius algebras}

The bilinear form on a Frobenius algebra endows it with a natural topology, turning it into a locally convex topological vector space. As we show in \Cref{lem:antipode-cont}, given a Hopf-Frobenius module, we can pull back this topology onto the Hopf algebra. This is what will allow us to give the ``Lie group'' associated with a Hopf-Frobenius module a topological group structure.

\newthought{Convention: } For the remainder of this section we will assume our base field $k$ is $\mathbb{R}$ or $\mathbb{C}.$

\begin{lemma}{Frobenius Algebra Topology is Induced from Vector Space Basis}{frob-top-basis}
Let $F$ be a vector space with bilinear form $\gen{-,-}$, and $B = \{ b_i \}$ be any basis and let 
\begin{align*}
X_F &= \{ \gen{-,x}, \gen{x,-} \mid x \in F \} \\
X_B &= \{ \gen{-, b_i}, \gen{b_i, -} \mid b_i \in B \}.
\end{align*}
If $\tau_F$ is the initial topology on $F$ induced by the maps in $X_F$ and $\tau_B$ is the initial topology induced by $X_B$, then
\begin{align*}
\tau_B = \tau_F.
\end{align*}
\end{lemma}
\begin{proof-lemma}
Since $X_B \subseteq X_F$, $\tau_B \subseteq \tau_F$. For any finite choice of elements $v_1, \ldots, v_t$, and $r > 0$, let
\begin{align*}
U_r(v_1, \ldots, v_t) \equiv \{ x \in F \mid \abs{\gen{x, v_i}}, \abs{\gen{v_i, x}} < r \text{ for all } i \}.
\end{align*}
Then $\{ U_r(v_1, \ldots, v_t) \mid r > 0, t \in \mathbb{N}, v_i \in F \}$ forms a neighbourhood basis of zero in $\tau_F$.

\newthought{$\tau_F \subseteq \tau_B$: } For any $U_r(v_1, \ldots, v_t)$ we can express each $v_i$ in the basis $B$. In particular, we will only need finitely many basis elements to do so, say $\{ b_{j_1}, \ldots b_{j_m} \}$. Write $v_i = \sum_{\ell} \alpha_{i\ell} b_{j_\ell}$. Now consider the open set
\begin{align*}
V = \{ x \mid \abs{\gen{x, b_{j_\ell}}}, \abs{\gen{b_{j_\ell}, x}} < \dfrac{r}{m (\max_{i\ell} \abs{\alpha_{i \ell}})}  \text{ for all } j_\ell \}.
\end{align*}
Then for any $x \in V$ we have
\begin{align*}
\abs{\gen{x, v_i}} &= \abs{\sum_{\ell} \alpha_{i \ell} \gen{x, b_{j_\ell}}} \\ 
&\leq \sum_{\ell} \abs{\alpha_{i \ell}} \abs{\gen{x, b_{j_\ell}}} \\
&< \sum_{\ell} \dfrac{r \abs{\alpha_{i \ell}} }{m (\max_{i\ell} \abs{\alpha_{i \ell}})} \\
&\leq r.
\end{align*}
Similarly, $\abs{\gen{v_i, x}} < r$ for all $i$. Thus $V \subseteq U_r(v_1, \ldots, v_t)$. It follows that $\tau_F \subseteq \tau_B$.
\end{proof-lemma}

\begin{lemma}{}{bilinear-cont-at-zero}
Let $V,W$ be topological vector spaces. Then a bilinear map $\gen{-,-}: V_1 \times V_2 \to W$ is continuous if and only if it is continuous at $(0,0)$.
\end{lemma}
\begin{proof-lemma}
Say $\gen{-,-}$ is continuous at $(0,0)$. Now take an arbitrary point $(x,y) \in V_1 \times V_2$, and an open neighbourhood $U$ of $\gen{x,y} \in W$. Since scalar multiplication and translation are homeomorphisms, $U - \gen{x,y}$ is an open neighbourhood of $0 \in W$, and there is a balanced open neighbourhood of $0$, $B \subseteq W$ such that
\begin{align*}
B + B + B \subseteq U - \gen{x,y}.
\end{align*}
Now since $\gen{-,-}$ is continuous at $(0,0)$ and $\gen{0,0} = 0$ by bilinearity, we can find a basic open set $A_1 \times A_2 \subseteq V_1 \times V_2$ such that
\begin{align*}
\gen{-,-} (A_1 \times A_2) \subseteq B.
\end{align*}
Without loss of generality, since $A_1, A_2$ are open sets containing zero, we can take them to be balanced and absorbing. Hence there are scalars 
\begin{align*}
0 < \lambda, \mu < 1
\end{align*}
such that $\lambda x \in A_1$ and $\mu y \in A_2$. Consider the open neighbourhood of $(x,y)$
\begin{align*}
\mu A_1 \times \lambda A_2 + (x,y).
\end{align*}
Given $(a,b) \in A_1 \times A_2$, we have
\begin{align*}
\gen{\mu a + x, \lambda b + y} &= \mu \lambda \gen{a,b} + \lambda \gen{x,b} + \mu \gen{a,y} + \gen{x,y} \\
&= \mu \lambda \gen{a,b} + \gen{\lambda x, b} + \gen{a, \mu y} + \gen{x,y}.
\end{align*}
Since $0 < \lambda, \mu < 1$ and $B$ is balanced, we have $\mu \lambda \gen{a,b} \subseteq \mu \lambda B \subseteq B$. We also have $\gen{\lambda x, b}, \gen{a, \mu y} \in B$. So
\begin{align*}
\gen{-,-}(\mu A_1 \times \lambda A_2 + (x,y)) &\subseteq B + B + B + \gen{x,y} \\
&\subseteq (U - \gen{x,y}) + \gen{x,y} \\
&= U.
\end{align*}
Thus $\gen{-,-}$ is continuous at $(x,y)$, and hence everywhere. The other direction of the proof is trivial.
\end{proof-lemma}

\begin{lemma}{Finite-Dimensional Frobenius Algebra Pairings are Jointly Continuous}{frob-jointly-cont}
Let $F$ be a finite-dimensional Frobenius algebra, viewed as a topological vector space with topology induced from $\gen{x,-}_{\text{Frob}}, \gen{-,y}_{\text{Frob}}$. Then the pairing
\begin{align*}
\gen{\cdot, \cdot}_{\text{Frob}}: F \times F \to k
\end{align*}
is jointly continuous. Similarly, if $F$ is over a trivially-starred field, and we give $F$ the topology induced by $\gen{x,-}, \gen{-,y}$, then
\begin{align*}
\gen{\cdot, \cdot}: F \times F \to k
\end{align*}
is jointly continuous. 
\end{lemma}
\begin{proof-lemma}
We prove that $\gen{\cdot, \cdot}_{\text{Frob}}$ is jointly continuous in the appropriate topology. The proof for $\gen{\cdot, \cdot}$ is analogous. 
Let $\{ b_i \}$ be a basis of $F$ and consider the following preimages of open balls around zero in $k$:
\begin{align*}
B^i_R(r) &\equiv \gen{-, b_i}_{\text{Frob}}^{-1}B_0(r)\\
B^i_L(r) &\equiv \gen{b_i,-}_{\text{Frob}}^{-1}B_0(r).
\end{align*}

\newthought{Basis: } From \Cref{lem:frob-top-basis} a neighbourhood basis of zero in $F$ is given by $\{ \bigcap_{i=1}^n B_L^i(r) \cap B_R^i(r) \mid $ for all $r > 0 \}$.

\newthought{Joint continuity: } From \Cref{lem:bilinear-cont-at-zero}, it's enough to prove continuity at $(0,0) \in F \times F$. Consider a ball of radius $r$ around zero $B_0(r)$ in $k$. We take the open neighbourhood of $(0,0)$
\begin{align*}
V = (\bigcap_{i=1}^n B_L^i(r) \cap B_R^i(r)) \times (\bigcap_{i=1}^n B_L^i(r) \cap B_R^i(r)).
\end{align*}
Then $\gen{\cdot, \cdot}_{\text{Frob}}(V) \subseteq B_0(r)$. Thus the pairing $\gen{\cdot, \cdot}_{\text{Frob}}$ is continuous at $(0,0)$ and hence is jointly continuous.
\end{proof-lemma}

\section{Examples of Hopf-Frobenius modules}
\label{sec:hf-mod-exs}

Part of the motivation behind the definition of Hopf-Frobenius modules is their ability to unify constructions that appear in geometry, functional analysis, and algebra. Though in the present paper we utilize only the finite-dimensional structure in our version of the Lie correspondence (and hence don't explore the functional-analytic side in detail), we make note of analytic examples of Hopf-Frobenius modules below, with the intention to utilize these modules in future work studying the case of infinite-dimensional Lie groups.

\subsection{Examples from geometry and algebra}

The intuition behind the examples relevant to Lie theory is that the Hopf algebra contains the data of the Lie group or Lie algebra, while the Frobenius algebra corresponds to the underlying manifold being acted on. The Hopf-Frobenius module conditions appear as the condition that the actions of the Lie group or Lie algebra must behave like symmetries or derivations, respectively.

\begin{lemma}{Hopf-Frobenius Modules Can be Pulled Back}{}
Let $(H,F)$ be a Hopf-Frobenius module, let $G$ be a (graded) Hopf algebra, and let $\phi: G \rightarrow H$ be a (graded) Hopf algebra homomorphism (of degree zero). Then $(G, \phi^{*} F)$ is a Hopf-Frobenius module, where $\phi^{*} F = F$ is restriction of scalars on the module $F$, viewed as a $G$-module:
\begin{equation*}
 g \cdot x = \phi(g) \cdot x,
\end{equation*}
for all $g \in G$ and $x \in F$.

\end{lemma}
\begin{proof-lemma}
By standard module theory,  $(G, \phi^{*} F)$ is a left $G$-module. We now check the other two Hopf-Frobenius conditions. Using the fact that $\phi$ is degree-preserving, we have
  \begin{align*}
  \sum (-1)^{\abs{x} \abs{g_{(2)}}} \gen{g_{(1)}\cdot x, g_{(2)}\cdot y} &= \sum (-1)^{\abs{x} \abs{\phi(g_{(2)})}} \gen{\phi(g_{(1)}) x, \phi(g_{(2)}) y}\\ 
  &=\sum  (-1)^{\abs{x} \abs{\phi(g_{(2)})}} \gen{\phi(g)_{(1)} x, \phi(g)_{(2)} y}\\
  &= \varepsilon_{H}(\phi(g)) \gen{x,y}\\
  &= \varepsilon_{G}(g) \gen{x,y},
  \end{align*}
  where in the second and final lines we used the fact that $\phi$ is a morphism of Hopf algebras, and thus commutes with the coproduct and counit, respectively. Similarly, we have
  \begin{align*}
   \gen{g \cdot x, y} &= \gen{\phi(g) x, y} = (-1)^{\abs{\phi(g)} \abs{x}} \gen{x, S_{H}(\phi(g))y}\\ 
 &=   (-1)^{\abs{g} \abs{x}} \gen{x, \phi(S_{G}(g))y}\\
 &=   (-1)^{\abs{g} \abs{x}} \gen{x, S_{G}(g) \cdot  y}
  \end{align*}
and so $(G, \phi^{*} F)$ is a Hopf-Frobenius module.
\end{proof-lemma}

\begin{example}{HF-Module: Trivial Hopf Action}{}
  Let $F$ be any Frobenius algebra and let $H$ be any Hopf algebra, both over the field $k$. Then there is a trivial Hopf-Frobenius module structure on $(H,F)$ given by pulling back the counit
\begin{equation*}
\varepsilon : H \rightarrow k.
\end{equation*}
Explicitly, $h \cdot f = \varepsilon(h) f$ for all $h \in H$, $f \in F$. 
\end{example}

\begin{example}{HF-Module: Lie Groups}{}
Let $O(n)$ be the orthogonal group, and take $F = \mathbb{R}^n$ with the standard basis and let $H = \mathbb{R} O(n)$ be the group algebra of the orthogonal group. Letting $O(n)$ act on $F$ by its standard representation and extending to $H$ linearly, gives an HF-module.
\end{example}
\begin{exampleframe}
We're grouplike-generated, so the HF-conditions become
\begin{align*}
\gen{g x,g y} = \gen{x,y},
\end{align*}
which is just the condition that $O(n)$ acts orthogonally on $\mathbb{R}^n$.

\newthought{Similarly: } We can construct HF-modules corresponding to most of the classical Lie groups: $SO(n)$, $U(n)$, $SU(n)$, $O(m,n)$, $SO(m,n)$, and the symplectic groups.

\newthought{Symplectic: } We'll consider the case of the symplectic groups in more detail, since there are a number of variations of its definition, depending on the base field. The real symplectic group $Sp(\mathbb{R}^{2n}, \mathbb{R})$ is defined to be group of invertible linear transformations on $\mathbb{R}^{2n}$ which preserve the symplectic pairing of \Cref{ex:frb-symp}. If you replace $\mathbb{R}^{2n}$ by $\mathbb{C}^{2n}$ in \Cref{ex:frb-symp}, the maps preserving the symplectic pairing form the complex symplectic group $Sp(\mathbb{C}^{2n}, \mathbb{C})$. Confusingly, the quaternionic symplectic group is not obtained by replacing $\mathbb{R}^{2n}$ by $\mathbb{H}^{2n}$, as the non-commutativity of the base (skew) field $\mathbb{H}$ prevents the symplectic pairing from being bilinear (and replacing the base field with $\mathbb{R}$ or $\mathbb{C}$ will just recover the real and complex symplectic groups, respectively). Instead, the quaternionic symplectic group is defined to be the group of transformations on $\mathbb{H}^{n}$ (not $2n$) preserving the quaternionic pairing
\begin{equation*}
 \gen{x,y} = \sum_{i} x_{i}^{*} y_{i}, \quad x, y \in \mathbb{H}^{n}.
\end{equation*}
This symplectic group is usually denoted $Sp(n)$. We note that $\mathbb{H}^{n}$ is also a Frobenius algebra by the generalization at the bottom of \Cref{ex:frb-euclid}.

\newthought{Note: } The above Lie groups fit naturally into the Hopf-Frobenius framework because they preserve certain bilinear forms. In the next section, we will show how, for these types of form-preserving Lie groups, some of the main structure of Lie theory can be obtained simply by working within the Hopf-Frobenius world, without appealing to any explicit differential geometry.
\end{exampleframe}

\begin{example}{HF-Module: Lie Algebras}{}
Let $\mathfrak{o}(n) = \mathfrak{so}(n)$ be the orthogonal Lie algebra (i.e. the skew-symmetric matrices) and take our Hopf algebra to be the universal enveloping algebra $U[\mathfrak{o}(n)]$. Let $F = \mathbb{R}^n$ with the standard basis. Let $\mathfrak{o}(n)$ act on $F$ by its standard representation and extend to $H$ by letting products of elements of $\mathfrak{o}(n)$ act on $F$ by composition. This has the structure of an HF-module.
\end{example}
\begin{exampleframe}
Being primitively generated, the HF-conditions are
\begin{align*}
\gen{\phi x,y} = -\gen{x,\phi y}.
\end{align*}
for $\phi \in \mathfrak{o}(n)$. This says that elements of $\mathfrak{o}(n)$ are represented by skew-symmetric matrices, which is precisely what happens in the standard representation!

\newthought{Similarly: } We can construct HF-modules corresponding to Lie algebras for the classical groups mentioned in the previous example.
\end{exampleframe}

\subsection{Examples from functional analysis}

In the infinite-dimensional case, our Hopf algebra plays the same role, acting as symmetries or derivations, except in this setting, the Frobenius algebra is thought of as a function space.

\begin{example}{HF-Module: Schwartz Space with Derivatives}{}

Let $F = \cat{S}(\mathbb{R}^n)$ be Schwartz space as in \Cref{ex:frb-schwartz} and $H = \mathbb{R}[\partial_1, \ldots, \partial_n]$ be a polynomial algebra whose action on $F$ is given by partial derivatives. Then $(H,F)$ is a Hopf-Frobenius module.
\end{example}
\begin{exampleframe}
Since $H$ is primitively generated, the Hopf-Frobenius conditions reduce to
\begin{equation}
\gen{\partial_i f, g} = - \gen{f,\partial_i g},
\end{equation}
or in integral notation, $\int (\partial_i f) g = - \int f \partial_i g$. This is precisely the integration by parts formula for functions (like those in $F$) which go to zero at infinity.

\newthought{Similarly: } $(\mathbb{R}[\partial_1, \ldots, \partial_n], C_c^{\infty}(\mathbb{R}^n))$ is a Hopf-Frobenius module.
\end{exampleframe}

\begin{example}{HF-Module: Compactly Supported Functions with Derivations}{}
Let $M$ be an orientable smooth $n$-manifold equipped with a volume form $\omega \in \Omega^n(M)$. Take as Frobenius algebra $F = C^{\infty}_c(M)$ as in \Cref{ex:frb-comp-supp}. Consider the Lie algebra of derivations $\operatorname{Der}(C^{\infty}_c(M))$, which is isomorphic to the Lie algebra of vector fields $\mathfrak{X}(M)$, which acts on $C^{\infty}_c(M)$. Let $H = U[\operatorname{Der}(C^{\infty}_c(M))]$ be the universal enveloping algebra on the derivations so that $F$ is a left $H$-module. The pair $(H,F)$ is a Hopf-Frobenius module.
\end{example}

\begin{proposition}{}{}
Let $G$ be a locally compact Hausdorff topological group and consider the group Hopf algebra $H = kG$. Let $F = C_c(G)$ be the algebra of compactly supported functions on $G$, and give $F$ the structure of a left $H$-module by $C_c(G)$ by $g \cdot \phi(h) = \phi (hg)$. Then given any measure $\mu$ on $G$ such that the linear functional
\begin{equation}
\varepsilon_\mu(\phi) = \int \phi d\mu
\end{equation}
has a simple kernel, $(H,F)$ is a Hopf-Frobenius module if and only if $\mu$ is a Haar measure.
\end{proposition}
\begin{proof-prop}
Since we're working with a Hopf algebra generated by grouplike elements, the Hopf-Frobenius conditions reduce to
\begin{align*}
\gen{g \phi,g \psi} = \gen{\phi, \psi}.
\end{align*}
In integral notation, this is $\int \phi( xg) \psi( xg) d \mu(x) = \int \phi \psi$. Since we're working with compactly supported functions, given any such function, $\phi$, and a $g \in G$, take a function $1_\phi^g$ which is equal to one on at least the support of $\phi$ and the image of $\operatorname{supp} \phi$ under the action of $g^{-1}$, and goes to zero outside of them. Then
\begin{equation}
\int \phi(xg) d \mu(x) = \int \phi(xg) 1_\phi^g (xg) d \mu(x) = \int \phi 1_\phi^g = \int \phi d \mu.
\end{equation}
Thus $\varepsilon$ is left-invariant with respect to $G$. This is equivalent to $\mu$ being a Haar measure on $C_c(G)$. Conversely, if $\mu$ is a Haar measure, it is obvious that the condition $\gen{g \phi,g \psi} = \gen{\phi, \psi}$ will be satisfied.
\end{proof-prop}

%

\section{The Lie Group-Lie Algebra Correspondence for Hopf-Frobenius Modules}
\label{sec:lie-correspondence}

In this section we recover the Lie group-Lie algebra correspondence for a certain class of Lie groups and algebras using the structure of Hopf-Frobenius modules. This allows us to motivate a definition of a topological group and Lie algebra associated to any Hopf-Frobenius module.

\newthought{Convention: } From this point, we consider only the ungraded setting.

\subsection{The Lie algebra of a Hopf-Frobenius module}

The Lie algebra that we associate to a group will have no direct knowledge of its topological structure --- it only sees $G$ as a discrete group.

\begin{lemma}{}{lie-alg-for-grouplike}
Let $G$ be a group and consider the group algebra $kG$. The vector subspace of $kG$ generated by 
\begin{equation}
\{ g - g^{-1} \mid g \in G \}
\end{equation}
is a Lie algebra under the commutator bracket.
\end{lemma}
\begin{proof-lemma}
We check that the subspace is closed under the commutator. We have
\begin{align*}
[g - g^{-1}, h - h^{-1}] &= (g-g^{-1})(h - h^{-1}) - (h - h^{-1})(g-g^{-1}) \\
 &= gh - g^{-1} h - g h^{-1} + g^{-1} h^{-1} - (hg - h^{-1} g - h g^{-1} + h^{-1} g^{-1}) \\
  &= (gh  - h^{-1} g^{-1}) + (g^{-1} h^{-1} - hg) + (h^{-1} g - g^{-1} h) + (hg^{-1} - gh^{-1}).
\end{align*}
Each of the four bracketed terms is of the form $x - x^{-1}$ with $x \in G$.
\end{proof-lemma}

\medskip

 It is quite surprising that this construction will yield the appropriate Lie algebra for all of the Lie groups we consider. Before we show this, however, we note that this Lie algebra construction can be generalized to apply to \textit{any} Hopf algebra.

\begin{lemma}{Lie Algebra of a Hopf Algebra}{lie-alg-of-hopf}
Let $H$ be a Hopf algebra and consider the vector subspace generated by
\begin{equation}
\{ x - S(x) \mid x \in H \text{ and } S^{2}(x) = x \}.
\end{equation}
This is a Lie algebra under the commutator bracket.
\end{lemma}
\begin{proof-lemma}
  This is essentially the same proof as \Cref{lem:lie-alg-for-grouplike}, and follows from the fact that $S(xy) = S(y) S(x)$ for every element of a Hopf algebra. Note also that the set of elements on which $S$ is an involution is a vector subspace, by linearity of $S$: $S^{2}(\alpha x + \beta y) = \alpha S^{2}(x) + \beta S^{2}(y) = \alpha x + \beta y$. Applying the commutator, we have
\begin{align*}
[x - S(x), y - S(y)] &= [x,y] - [x, S(y)] - [S(x), y] + [S(x), S(y)] \\
&= [x,y] - [x, S(y)] - S([S(y), x]) + S([y,x]) \\
&= [x,y] - S([x,y]) - \left( [x,S(y)] - S([x, S(y)]) \right).
\end{align*}
Now the first pair of terms and the second pair are both of the form $z - S(z)$.
\end{proof-lemma}

\newthought{Note: } In $kG$ we have $S(g) = g^{-1}$, and since antipodes are linear, the linear combinations of all $\{x - S(x) \}$ will match the vector space generated by $\{ g - g^{-1} \}$ used in the previous construction.

\begin{definition}{Lie Algebra of a Hopf Algebra}{}
Let $H$ be a Hopf algebra. The \emph{Lie algebra associated to $H$} is the vector space
\begin{equation*}
\mathfrak{L}_H \equiv \{ \sum_i \alpha_i (x - S(x)) \mid x \in H, S^2(x) = x, \alpha_i \in k \},
\end{equation*}
with the commutator bracket.
\end{definition}
\begin{defframe}
When $G$ is a group and $H$ is the group Hopf algebra $kG$, we denote $\mathfrak{L}_{H}$ as $\mathfrak{L}_G$, if the field $k$ is clear.
\end{defframe}

\begin{proposition}{}{}
Let $\mathfrak{L}: \textsc{Gr} \to \textsc{L.Alg}$ be the assignment from the category of groups to the category of Lie algebras which sends $G$ to $\mathfrak{L}_G$ and $\phi: G \to H$ to the morphism $\mathfrak{L}_\phi$ defined by 
\begin{equation*}
\mathfrak{L}_\phi(g - g^{-1}) = \phi(g) - \phi(g)^{-1}.
\end{equation*}
This is a functor.
\end{proposition}
\begin{proof-prop}

\newthought{Well-defined: } Note that in $kG$ the grouplike elements are linearly independent \cite[Proposition~1.4.14]{dascalescutextbook}. It follows that the only linear dependence between different $g - g^{-1}$ is the relation $-(g - g^{-1}) = g^{-1} - (g^{-1})^{-1}$, which is respected by our definition of $\mathfrak{L}_\phi$. So $\mathfrak{L}_\phi$ is well-defined.

\newthought{Lie algebra map: } We now check that $\mathfrak{L}_{\phi}$ is a Lie algebra morphism. Using the fact that 
\begin{align*}
[g-g^{-1}, h-h^{-1}] = (gh  - h^{-1} g^{-1}) + (g^{-1} h^{-1} - hg) + (h^{-1} g - gh^{-1}) + (hg^{-1} - gh^{-1}),
\end{align*}
we have 
\begin{align*}
\mathfrak{L}_{\phi}([g-g^{-1}, h - h^{-1}]) &= \mathfrak{L}_\phi (gh  - h^{-1} g^{-1}) + \ldots + \mathfrak{L}_\phi (hg^{-1} - gh^{-1}) \\
&= (\phi(gh) - \phi(gh)^{-1}) + \ldots + (\phi(hg^{-1}) - \phi(gh^{-1})^{-1}) \\
&= (\phi(g) \phi(h) - \phi(h)^{-1} \phi(g)^{-1}) + \ldots + (\phi(h) \phi(g^{-1}) - \phi(g) \phi(h)^{-1}) \\
&= [\phi(g) - \phi(g)^{-1}, \phi(h) - \phi(h)^{-1}].
\end{align*}

\newthought{Functorial: } Let $G \xrightarrow{\phi} H \xrightarrow{\psi}K$. We have that $\mathfrak{L}_\psi \mathfrak{L}_\phi(g-g^{-1}) = \mathfrak{L}_\psi (\phi(g) - \phi(g)^{-1}) =  \psi (\phi(g)) - \psi(\phi(g))^{-1} = \mathfrak{L}_{\psi \circ \phi}(g - g^{-1})$.
\end{proof-prop}

\medskip

We next show that for Lie groups we're interested in, the above construction reproduces their standard Lie algebras.

\begin{proposition}{}{lie-alg-matches}
Let $G$ be a Lie group over a field $k$ that satisfies the following property: 

\begin{enumerate}
\item There exists a finite collection of non-degenerate bilinear forms $\gen{-,-}_{i}$ on a vector space $V$ such that $G = \{ g \in \operatorname{End}(V) \mid \gen{g x, gy}_i = \gen{x,y}_i$ for all $i \}$, and the Lie algebra of $G$ is $\mathfrak{g} = \{ A \in \operatorname{End}(V) \mid \gen{Ax, y}_i = - \gen{x, Ay}_i $ for all $i \}$.
\end{enumerate}

Then $\mathfrak{L}_G \cong \mathfrak{g}$.
\end{proposition}
\begin{proof-prop}
Given an element of the Lie group, $M$, and a bilinear form $\gen{\cdot, \cdot}_i$ we have 
\begin{align*}
\gen{(M - M^{-1})x, y}_i &= \gen{Mx, y}_i - \gen{M^{-1}x, y}_i \\
&= \gen{x, M^{-1} y}_i -  \gen{x, M y}_i \\
&= - \gen{x, (M - M^{-1}) y}_i.
\end{align*}
So we have an assignment
\begin{equation}
\mathfrak{L}: \mathfrak{L}_G \to \mathfrak{g}
\end{equation}
that sends $M - M^{-1} \mapsto M - M^{-1}$. This assignment induces a Lie algebra isomorphism $\mathfrak{L}_G \cong \mathfrak{g}$.
\end{proof-prop}

\begin{example}{Lie Groups Satisfying \Cref{prop:lie-alg-matches}}{}

This class of Lie groups and Lie algebras covers a wide range, including most of the classical Lie groups/algebras: $O(n)$, $U(n)$, $Sp(2n)$, $O(p,q)$ (the indefinite orthogonal group), and in particular, the Lorentz group $\Lambda = O(1,3)$. 

\end{example}
\begin{exampleframe}
\newthought{Other groups: } Several additional classical Lie groups which \textit{don't} satisfy the above condition nevertheless still reproduce the correct Lie algebra, as we outline in the following proposition.
\end{exampleframe}

\begin{proposition}{}{lie-alg-matches-gl-sl}
Let $G$ be $GL(n)$ or $SL(n)$ over a field $k$ with $\operatorname{char}(k) \ne 2$, with corresponding Lie algebra $\mathfrak{g}$. Then
\begin{equation*}
\mathfrak{L}_G \cong \mathfrak{g}.
\end{equation*}
\end{proposition}
\begin{proof-prop}

\newthought{$GL(n)$: } Consider the invertible matrices $M = I + E_{ij} \in GL(n)$, $i \ne j$. We have $M - M^{-1} = 2 E_{ij}$. Similarly, if $M = I + E_{ii} \in GL(n)$, then $M - M^{-1} = \dfrac{3}{2} E_{ii}$. Thus $\mathfrak{L}_G$ contains all of the generators of $\mathfrak{gl}(n)$, so $\mathfrak{L}_G \cong \mathfrak{gl}(n)$.

\newthought{$SL(n)$: } We again take $M = I + E_{ij} \in SL(n)$ to get the generators $2 E_{ij}$. Our other generators are obtained from $M = I - \dfrac{1}{2} E_{ii} +  E_{jj} \in SL(n)$. We have $M - M^{-1} = \dfrac{3}{2}(E_{jj} - E_{ii})$. This will produce all of the remaining generators (actually, only considering $M = I - \dfrac{1}{2} E_{11} +  E_{jj}$ will be enough to produce the generators), and thus we have $\mathfrak{L}_G \cong \mathfrak{sl}(n)$.
\end{proof-prop}

\subsection{The Lie group of a Hopf-Frobenius module}

We now turn to the Lie group associated to a Hopf-Frobenius module. Recalling \Cref{ex:frb-orthog} and \Cref{ex:frb-symp}, we will consider the situation in which $(H,F)$ forms a Hopf-Frobenius module under potentially more than one bilinear form on $F$. This will be useful in the case of the unitary group.

\begin{proposition}{The Group Generated by Primitive Elements}{lie-group-of-hopf}
Let $H$ be a Hopf algebra. Let $F$ be an algebra which is not necessarily unital nor associative, and which is a left $H$-module. And let $\{ \gen{\cdot, \cdot}_i \}_{i \in I}$ be finitely many bilinear forms which each give $F$ the structure of a Frobenius algebra and turn $(H,F)$ into a Hopf-Frobenius module (for the same fixed action by $H$). Let $P$ be the primitive elements of $H$. Let $U = U[P]$ be the subalgebra generated by $P$ and consider the induced action of $U[P]$ on $F$. Then
\begin{align*}
U[P]_F \equiv \{ g \in U/\operatorname{ann} F \mid &\gen{g x, gy}_i = \gen{S(g) x, S(g) y}_i = \gen{x,y}_i \\
 &\text{ and } \gen{S(g) x, y}_i = \gen{x, g y}_i \text{ for all } i \}
\end{align*}
forms a group, where $g^{-1} = S(g)$. Here  $\operatorname{ann} F \subseteq U$ denotes the annihilator of $F$ with respect to the $U$-module structure.
\end{proposition}
\begin{proof-prop}
First notice that $g \in U[P]_F$ if and only if $S(g) \in U[P]_F$.

\newthought{Closed under Products: } Given $g,h \in U[P]_F$, we have $\gen{gh x,gh y} = \gen{hx,hy} = \gen{x,y}$. Now we use the fact that the antipode is an antihomomorphism \cite[Proposition~4.2.6]{dascalescutextbook} to write
\begin{align*}
\gen{S(gh)x,S(gh)y} &= \gen{S(h) S(g),S(h)S(g)y} = \gen{x,y} \\
\gen{S(gh)x,y} &= \gen{S(h)S(g)x,y} = \gen{S(g)x,hy} = \gen{x,gh y}
\end{align*}
Thus $U[P]_F$ is closed under multiplication.

\newthought{Identity: } Let $1$ be the image of the unit of $H$ in $U/\operatorname{ann} F$. It clearly satisfies the conditions to be in $U[P]_F$ and is the unit of the group.

\newthought{Inverses: } Let $g \in U[P]_F$. By assumption $\gen{x,y} = \gen{g x, gy} = \gen{S(g) g x, y}$ for all $y$. Thus by non-degeneracy of the pairing, $(S(g) g - 1)x = 0$ for all $x$. Since we've quotiented out the annihilator, this implies $S(g) g = 1$. Similarly, $\gen{S(g) x, S(g) y} = \gen{x, g S(g) y}$, giving $g S(g) = 1$. Thus
\begin{equation*}
g^{-1} = S(g). \qedhere
\end{equation*}
\end{proof-prop}

\newthought{Convention: } We will denote Hopf-Frobenius modules of the form that appears in \Cref{prop:lie-group-of-hopf} by 
\begin{align*}
(H,F, \{ \gen{\cdot, \cdot}_i \}_{i \in I}).
\end{align*}
In the case of only a single Frobenius algebra structure on $F$, this reduces to the ordinary notion.

\begin{definition}{Lie Group of a Hopf-Frobenius Module}{lie-group-hf-mod}
Let $(H,F, \{ \gen{\cdot, \cdot}_i \}_{i \in I})$ be as above and let $U$ be the subalgebra generated by the primitive elements of $H$. We denote the group constructed above by
\begin{align*}
\cat{G}^H_F \equiv \{ g \in U/\operatorname{ann} F \mid &\gen{g x, gy}_i = \gen{S(g) x, S(g) y} = \gen{x,y}_i \\
 &\text{ and } \gen{S(g) x, y}_i = \gen{x, g y}_i, \text{ for all } i \}
\end{align*}
and call it the \emph{Lie group of $(H,F, \{ \gen{\cdot, \cdot}_i \}_{i \in I})$}.
\end{definition}
\begin{defframe}

  \newthought{Smooth structure: } A Lie group should be more than a discrete group. By \cite[Corollary~3.41]{warner1971} we know that topological groups have at most one Lie group structure up to isomorphism. As such, we aim to view the Lie group of a Hopf-Frobenius module just as a topological group. We do such in the next sequence of proofs.
\end{defframe}

\begin{lemma}{}{lie-group-gives-hfm}
Let $(H,F, \{ \gen{\cdot, \cdot}_i \}_{i \in I})$ be a Hopf-Frobenius module. Then $(k \cat{G}_F^H, F, \{ \gen{\cdot, \cdot}_i \}_{i \in I})$ is also a Hopf-Frobenius module.
\end{lemma}
\begin{proof-lemma}
Since $U[P]$ acts on $F$ and $\operatorname{ann} F$ is an ideal of $U[P]$, $F$ is a left $\cat{G}_F^H$-module. The Hopf-Frobenius conditions are satisfied by definition, since for all $g \in \cat{G}_F^H$ and $i \in I$
\begin{align*}
\gen{gx, gy}_i &= \gen{x,y}_i \\
\gen{S(g)x, y}_i &= \gen{x, gy}_i. \qedhere
\end{align*}
\end{proof-lemma}

\begin{lemma}{}{antipode-cont}
Let $(H,F, \{ \gen{\cdot, \cdot}_i \}_{i \in I})$ be a Hopf-Frobenius module with Lie group $\cat{G}^H_F$, where $F$ is given the topology induced by all $\gen{x,-}_i, \gen{-,y}_i$. Give $\cat{G}^H_F$ the initial topology induced by the maps
\begin{equation}
\cat{G}^H_F \xrightarrow{- \cdot x} F
\end{equation}
for all $x \in F$. Then the antipode map $S: \cat{G}^H_F \to \cat{G}^H_F$ is continuous.
\end{lemma}
\begin{proof-lemma}

\newthought{Antipode: } From \Cref{prop:lie-group-of-hopf}, the antipode map is the inverse map on $\cat{G}^H_F$, and in particular is involutive. Now since the topology on $\cat{G}^H_F$ is initial, we need to check that the compositions $\cat{G}_F^H \xrightarrow{S} \cat{G}_F^H \xrightarrow{- \cdot x} F$ are continuous. These are continuous if and only if the postcompositions with $\gen{a,-}_i$ and $\gen{-,b}_i$ are continuous for all $a,b \in F$. We have $\gen{a, S(g) x}_i = \gen{S^2(g)a, x}_i = \gen{g a, x}_i$ and $\gen{S(g)x, b}_i = \gen{x, gb}_i$. So we're looking at the following commutative diagram:
\begin{center}
\begin{tikzcd}
& F \ar[r, "{\gen{\cdot, x}}_i"] & k \\
\cat{G}_F^H \ar[r, "S"] \ar[dr, swap, "{- \cdot x}"] \ar[ur, "{- \cdot a}" ] & \cat{G}_F^H  \ar[r, "{- \cdot x}"] & F \ar[d, "{\gen{\cdot, b}}_i"] \ar[u, swap, "{\gen{a,\cdot}}_i"]  \\
& F \ar[r, "{\gen{x,\cdot}}_i"] & k
\end{tikzcd}
\end{center}
The outer maps are continuous by construction of the topologies on $F$ and $\cat{G}_F^H$. Thus $S$ is continuous.
\end{proof-lemma}

\begin{proposition}{}{action-cont}
Let $(H,F, \{ \gen{\cdot, \cdot}_i \}_{i \in I})$ be a Hopf-Frobenius module where $F$ is finite-dimensional and over a trivially-starred field. Then the action of $\cat{G}_F^H$ on $F$
\begin{align*}
- \cdot - : \cat{G}_F^H \times F \to F
\end{align*}
is jointly continuous.
\end{proposition}
\begin{proof-prop}
By the characteristic property of the initial topology, the action is continuous if and only if its post composition with each $\gen{x,-}_i$ and $\gen{-, x}_i$ is continuous. We have
\begin{align*}
\gen{x, -}_i \circ (- \cdot -)(g,y) = \gen{x, g \cdot y}_i = \gen{S(g) \cdot x, y}_i.
\end{align*}
Thus the following diagram commutes
\begin{center}
\begin{tikzcd}
\cat{G}_F^H \times F \ar[r, "S \times 1"] \ar[d, swap, "{(- \cdot  -)}"] & \cat{G}_F^H \times F  \ar[r, "{(- \cdot x) \times 1}"]& F \times F \ar[d, "{\gen{\cdot, \cdot}}_i"] \\
F \ar[rr, swap, "{\gen{x, -}_{i}}"] && F
\end{tikzcd}
\end{center}
By the definition of the topology on $\cat{G}_F^H$, \Cref{lem:frob-jointly-cont}, and \Cref{lem:antipode-cont} respectively, each of $(- \cdot x)$,~$\gen{\cdot, \cdot}_i$, and $S$, are continuous, and thus $\gen{x,-}_i \circ (- \cdot -)$ is continuous as well. Similarly, $\gen{g \cdot y, x}_i = \gen{y, S(g) \cdot x}_i$, so the diagram
\begin{center}
\begin{tikzcd}
\cat{G}_F^H \times F \ar[r, "S \times 1"] \ar[d, swap, "{(- \cdot -)}"] & \cat{G}_F^H \times F  \ar[r, "{(\cdot, x) \times 1}"]& F \times F  \ar[r, "\tau"] & F \times F \ar[d, "{\gen{\cdot, \cdot}}_i"] \\
F \ar[rrr, swap, "{\gen{-,x}}_i"] &&& F 
\end{tikzcd}
\end{center}
commutes, showing that $\gen{-,x}_i  \circ (- \cdot -)$ is continuous. It follows that the action map is jointly continuous.
\end{proof-prop}

\begin{proposition}{The Lie Group of an HF-Module is a Topological Group}{top-group}
Let $(H,F, \{ \gen{\cdot, \cdot}_i \}_{i \in I})$ be a Hopf-Frobenius module where $F$ is finite-dimensional and over a trivially-starred field. The initial topology induced by the maps
\begin{equation}
\cat{G}^H_F \xrightarrow{- \cdot x} F
\end{equation}
for all $x \in F$ gives $\cat{G}^H_F$ the structure of a topological group.
\end{proposition}
\begin{proof-prop}

\newthought{Inverse: } From \Cref{lem:antipode-cont}, the inverse (i.e. antipode) is a continuous map.

\newthought{Multiplication: } By the characteristic property of the weak topology, $\mu: \cat{G}_F^H \times \cat{G}_F^H \to \cat{G}_F^H$ is continuous if and only if $(- \cdot x) \circ \mu: \cat{G}_F^H \times \cat{G}_F^H \to F$ is continuous for all $x \in F$. Consider the diagram
\begin{center}
\begin{tikzcd}
\cat{G}_F^H \times \cat{G}_F^H  \ar[r, "\mu"] \ar[d, swap, "1 \times {(- \cdot x)}"] & \cat{G}_F^H \ar[d, "{(- \cdot x)}"] \\
\cat{G}_F^H \times F \ar[r, swap, "{(- \cdot -)}"] & F
\end{tikzcd}
\end{center}
This commutes due to the associativity of the Hopf-Frobenius action: $(g,h) \to (gh) \cdot x = g \cdot (h \cdot x)$. The left map of the above diagram is continuous by definition, and the bottom map is continuous from \Cref{prop:action-cont}. It follows that the multiplication $\mu$ is continuous.
\end{proof-prop}

\subsection{The Hopf-Frobenius-Lie correspondence}

We now begin to prove our main theorem: that the Lie algebras and topological groups we defined in the previous sections satisfy a kind of Lie correspondence for any Hopf-Frobenius module.

\begin{lemma}{}{filtration-gives-endo}
Let $\rho: \cat{A} \to \operatorname{End}(V)$ be a finite-dimensional faithful irreducible representation of a unital associative algebra over a field $k$ such that there is a filtration of the form
\begin{equation*}
\rho(A_m \ldots A_1)(V) \subseteq \rho(A_{m-1} \ldots A_1)(V) \subseteq \ldots \subseteq \rho(A_1)(V) \subseteq V
\end{equation*}
where $A_i \in \cat{A}$ and $\dim \rho(A_m \ldots A_1)(V) = 1$. Then
\begin{equation}
\rho(\cat{A}) \cong \operatorname{End}(V).
\end{equation}
\end{lemma}
\begin{proof-lemma}

\newthought{Transitivity: } Note that since our representation is irreducible, the action of $\cat{A}$ is transitive on $V$, since given any non-zero $v \in V$, the subspace
\begin{align*}
\{ A v \mid A \in \cat{A} \}
\end{align*}
is invariant under the action of $\cat{A}$, and thus must be equal to all of $V$. Since we can produce a rank one matrix, $\rho(A_m \ldots A_1)$, it follows that all rank one matrices can be obtained in the image of $\rho$. Since all matrices are in the span of rank one matrices, it follows that the map $\rho$ is bijective, and thus $\cat{A} \cong \operatorname{End}(V).$
\end{proof-lemma}

\begin{lemma}{}{lie-groups-that-work}
Let $G$ be one of the following Lie groups with its standard representation:
\begin{align*}
\{ O(n), SO(n) \mid n > 2 \} \cup \{ Sp(n), O(m,n) \}.
\end{align*}
Then the induced representation of the universal enveloping algebra $U[\mathfrak{g}]$ surjects onto $\operatorname{End}(V)$ for the corresponding $V$.
\end{lemma}
\begin{proof-lemma}
For each case, by \Cref{lem:filtration-gives-endo}, it's sufficient to find an appropriate filtration by images of elements of $U[\mathfrak{g}]$. Note that each of the standard representations of $U[\mathfrak{g}]$ is irreducible. In each case that follows we fix some basis of our vector space and view the representations of $U[\mathfrak{g}]$ as matrix representations.

\newthought{$O(n)$: } For $n > 2$, the Lie algebra $\mathfrak{o}(n)$ in its standard representation has elements of the form $E_{ij} - E_{ji}$, and $E_{jk} - E_{kj}$, with $i,j,k$ distinct. The product 
\begin{align*}
(E_{ij} - E_{ji})(E_{jk} - E_{kj}) = -E_{jk}
\end{align*}
is rank one. The same example holds for $\mathfrak{so}(n)$, since it's isomorphic to $\mathfrak{o}(n)$.

\newthought{$Sp(\mathbb{R}^{2n}, \mathbb{R})$: } The Lie algebra consists of block matrices of the form
\begin{align*}
\begin{pmatrix}
A & B \\
C & -A^T
\end{pmatrix},
\end{align*}
where $B$ and $C$ are symmetric matrices. In particular, we can take $B = E_{ii}$ and every other block to be zero, to get a rank one matrix.

\newthought{$O(m,n)$: }
\begin{align*}
\begin{pmatrix}
A & B \\
B^{T} & D
\end{pmatrix},
\end{align*}
where $A$ is $m \times m$, $D$ is $n \times n$, and both are skew-symmetric, and $B$ is an arbitrary real matrix. We can obtain any $E_{ij}$ in the $A$ and $D$ blocks using the argument for $\mathfrak{so}(n)$ above. Since $B$ is arbitrary, we can acquire all of the elements $E_{ij} + E_{ji}$ for $1 \leq i \leq m$ and $m < j \leq m + n.$ Given such $i$ and $j$, $E_{ii}$ will be in the $A$ block, and $E_{jj}$ in the $D$, block, allowing us to acquire the individual entries in the $B$ blocks, as follows:
\begin{align*}
  E_{ii}(E_{ij} + E_{ji}) &= E_{ij} \\
  E_{jj}(E_{ij} + E_{ji}) &= E_{ji}.
                            \qedhere
\end{align*}
\end{proof-lemma}

\medskip

We consider the unitary case separately, as there is a subtlety involved in looking at real versus complex spaces:

\begin{lemma}{}{unitary-surj}
Let $G = U(n)$ with $n \geq 2$ or $SU(n)$ with $n > 2$ and consider the Hopf-Frobenius module given by $\mathbb{R} G$ acting on $\mathbb{R}^{2n}$. By abuse of notation, let $\operatorname{End}(\mathbb{C}^{n})$ be the subspace of $\operatorname{End}(\mathbb{R}^{2n})$ identified with complex matrices.Then $U[\mathfrak{u}(n)]$ surjects onto $\operatorname{End}(\mathbb{C}^{n})$.
\end{lemma}
\begin{proof-lemma}
  For ease of notation, we identify the real matrices in the representation of $\mathfrak{u}(n)$ with their isomorphic complex matrices.

\newthought{$U(n)$: } The Lie algebra $\mathfrak{u}(n)$ consists of skew-Hermitian matrices. If $n \geq 2$, we can obtain all off-diagonal entries by using
\begin{align*}
(i E_{aa}) (E_{ab} - E_{ba}) &= i E_{ab} \\
(i E_{aa})(i E_{ab} + i E_{ba}) &= - E_{ab}.
\end{align*}
All diagonal entries are obtained from $i E_{aa}$ and
\begin{equation*}
(i E_{aa})^2 = - E_{aa}.
\end{equation*}

\newthought{$SU(n)$: } For the case of $\mathfrak{su}(n)$ we use
\begin{align*}
(i E_{aa} - \dfrac{i}{n} I) (E_{ab} - E_{ba}) &= (i - \dfrac{i}{n}) E_{ab} + \dfrac{i}{n} E_{ba} = \dfrac{(n-1)}{n} i E_{ab} + \dfrac{1}{n} i E_{ba} \\
(i E_{aa} - \dfrac{i}{n}I)(i E_{ab} + i E_{ba}) &= (-1 + \dfrac{1}{n}) E_{ab} + \dfrac{1}{n}E_{ba}.
\end{align*}
As long as $n  > 2$, we can write
\begin{align*}
  ( \dfrac{(n-1)}{n} i E_{ab} + \dfrac{1}{n} i E_{ba}) - \dfrac{1}{n} (iE_{ab} + iE_{ba}) &= (\dfrac{n-2}{n}) i E_{ab} \\
((-1 + \dfrac{1}{n}) E_{ab} + \dfrac{1}{n}E_{ba}) + \dfrac{1}{n}(E_{ab} - E_{ba}) &= (\dfrac{2 - n}{n}) E_{ab}
\end{align*}
to recover all of the off-diagonal elements $E_{ab}$ and $i E_{ab}$. The diagonal elements are then obtained from
\begin{align*}
  E_{ab}(E_{ab} - E_{ba}) &= - E_{aa} \\
  E_{ab}(i E_{ab} + i E_{ba}) &= i E_{aa}.
                                \qedhere
\end{align*}
\end{proof-lemma}

\begin{theorem}{}{lie-group-matches-topologically}
Let $G$ be one of the following Lie groups: $Sp(\mathbb{R}^{2n}, \mathbb{R})$, $\{ U(n) \mid n \geq 2 \}$, or $\{ O(n) \mid n > 2 \}$, and let $\rho: G \to GL(F)$ be its standard representation on the corresponding vector space $F$ with bilinear forms $\{ \gen{\cdot, \cdot}_i \}_{i \in I}$. In particular, the action of $G$ preserves the forms:
\begin{equation*}
\gen{gx,gy}_i = \gen{x,y}_i, \text{ for all } g \in G, x,y \in F, i \in I.
\end{equation*}
If $\mathfrak{g}$ is the Lie algebra of $G$ and $U[\mathfrak{g}]$ is its universal enveloping algebra then
\begin{equation*}
\cat{G}_F^{U[\mathfrak{g}]} \cong G
\end{equation*}
as topological groups.
\end{theorem}
\begin{proof-theorem}
First note that in all of the cases in \Cref{lem:lie-groups-that-work} our Frobenius algebra is finite-dimensional and over a trivially-starred field.

\newthought{Hopf-Frobenius module: } The representation $\rho: G \to GL(F)$ induces a Lie algebra representation
\begin{equation}
d\rho: \mathfrak{g} \to \mathfrak{gl}(F).
\end{equation}
It is from this representation that we get the structure of a Hopf-Frobenius module from the universal enveloping algebra of $\mathfrak{g}$ acting on $F$. Note that the primitive elements of the universal enveloping algebra $U[\mathfrak{g}]$ are precisely the elements of $\mathfrak{g}$, and these generate $U[\mathfrak{g}]$ as an algebra (taking all products and sums). So in this case $U[P]$, where $P$ is the set of primitive elements of the Hopf algebra, which we use to construct $\cat{G}_F^{U[\mathfrak{g}]}$, is the universal enveloping algebra of $\mathfrak{g}$. From \Cref{lem:lie-groups-that-work} we have a surjection
\begin{align*}
U[\mathfrak{g}] \to \operatorname{Mat}_n(k)
\end{align*}
for the cases of $O(n)$ and $SO(n)$ (with $n > 2$) and $Sp(n)$. Since $F$ is acted upon by $\cat{G}_F^{U[\mathfrak{g}]}$, and each element of $\cat{G}_F^{U[\mathfrak{g}]}$ preserves the pairing, we have a group homomorphism $\phi: \cat{G}_F^{U[\mathfrak{g}]} \to G$. Now given any $g \in G$, we can view it as its matrix representation $\rho(g)$, which will correspond to some element in $\operatorname{Mat}_n(k)$ which preserves the pairing. It follows from the above surjection that there is some $\phi_g \in \cat{G}_F^{U[\mathfrak{g}]}$, which is mapped to $g$, and thus our group homomorphism is an isomorphism.

\newthought{Unitary case: } In the case of $U(n)$ and $SU(n)$, we don't have a surjection from $U[\mathfrak{u}(n)]$ onto all of the real matrices, but by \Cref{lem:unitary-surj} we still acquire all \textit{unitary} matrices.

\newthought{Continuity: } The standard representation of the group $G$ is continuous, where the topology on matrices is the standard subspace topology induced from $k^{n^2}$. The topology on $\cat{G}_F^{U[\mathfrak{g}]}$ is homeomorphic to the subspace topology of the matrix topology. It follows that our bijection
\begin{align*}
G \to \cat{G}_F^{U[\mathfrak{g}]}
\end{align*}
is an isomorphism of topological groups.
\end{proof-theorem}

\medskip

\newthought{Special groups: } Note that if we apply the above construction to $SU(n)$ or $SO(n)$, $n > 2$, the Lie group $\cat{G}_F^{U[\mathfrak{g}]}$ will recover $U(n)$ and $O(n)$, respectively, as the group $\mathcal{G}$ is defined to pick out all of the elements that preserve the pairing, without any kind of determinant condition. In this way, the Hopf-Frobenius module version of the Lie correspondence retains full knowledge of the inner product space, recovering the entirety of the original isometry Lie group.

\subsection{Exceptional Lie Groups}

One should be able to adapt the construction of the group in \Cref{def:lie-group-hf-mod} to recover the five exceptional Lie groups from the universal enveloping algebras of their  corresponding Lie algebras, by noting that each of the exceptional Lie groups can be defined in terms of the octonions, which are Frobenius algebras by \Cref{ex:comp-alg}
